\documentclass[MSNbibl,nameyear,dvips]{arxbj}

\aid{0}
\volume{19}
\issue{4}
\pubyear{2013}
\firstpage{1088}
\lastpage{1121}
\doi{10.3150/12-BEJSP12} 

\makeatletter

\renewcommand{\Sigma}{\sum}

\renewcommand{\epsilon}{\varepsilon}

\makeatother

\begin{document}
\begin{frontmatter}

\title{A Tricentenary history of the Law of Large Numbers}
\runtitle{Tricentenary history}

\begin{aug}
\author{\fnms{Eugene} \snm{Seneta}\corref{}\ead[label=e1]{eseneta@maths.usyd.edu.au}}
\runauthor{E. Seneta} 
\address{School of Mathematics and Statistics FO7, University of Sydney, NSW 2006, Australia.\\
\printead{e1}}
\end{aug}


%
\begin{abstract}
The Weak Law of Large Numbers is traced chronologically from its inception
as Jacob Bernoulli's Theorem in 1713, through De Moivre's Theorem, to ultimate
forms due to Uspensky and Khinchin in the 1930s, and beyond. Both
aspects of Jacob Bernoulli's Theorem: 1. As limit theorem (sample size $n \to\infty$),
and: 2. Determining sufficiently large sample size for
specified precision, for known and also unknown $p$ (the inversion problem),
are studied, in frequentist and Bayesian settings. The Bienaym\'e--Chebyshev
Inequality is shown to be a meeting point of the French and Russian directions
in the history. Particular emphasis is given to less well-known aspects
especially of the Russian direction, with the work of Chebyshev, Markov (the
organizer of Bicentennial celebrations), and S.N. Bernstein as focal points.
\end{abstract}



\pdfkeywords{Bienayme--Chebyshev Inequality,
Jacob Bernoulli's Theorem,
J.V. Uspensky and S.N. Bernstein,
Markov's Theorem,
P.A. Nekrasov and A.A. Markov,
Stirling's approximation}

\end{frontmatter}

\section{Introduction}
\subsection{Jacob Bernoulli's Theorem}

Jacob Bernoulli's Theorem was much more than the first instance of what
came to be know in later times as the Weak Law of Large Numbers (WLLN).
In modern notation Bernoulli showed that, for fixed $p$, any given
small positive number $\epsilon$, and any given large positive number
$c$ (for example $c =1000$), $n$ may be specified so that:
%
\begin{equation}
\label{bernoulli1} P\biggl(\biggl| \frac{X}{n} - p\biggr| > \epsilon\biggr) <
\frac{1}{c+1}
\end{equation}
for $n \geq n_0(\epsilon, c)$. The context: $X$ is the number of
successes in $n$ binomial trials relating to sampling with replacement
from a collection of $r+s$ items, of which $r$ were ``fertile'' and $s$
``sterile'', so that $p= r/(r+s)$. $ \epsilon$ was taken as $1/(r+s)$.
His conclusion was that $n_0(\epsilon, c)$ could be taken as the
integer greater than or equal
to:
%
\begin{eqnarray}
\label{bernoulli2}
&&t \max {\biggl\{} \frac{\log c(s-1)}{\log(r+1) - \log r} {\biggl(} 1 +
\frac{s}{r+1} {\biggr)} - \frac{s}{r+1},\nonumber\\[-8pt]\\[-8pt]
&&\hspace*{20.5pt}\quad \frac{\log c(r-1)}{\log
(s+1) - \log s} {\biggl(} 1 +
\frac{r}{s+1} {\biggr)} - \frac{r}{s+1} {\biggr\}}\nonumber
\end{eqnarray}
where $t =r+s$. The notation $c,r,s,t$ is Bernoulli's and the form of
the lower bound for $n_0(\epsilon, c)$ is largely his notation.

There is already a clear understanding of the concept of an event as a
subset of outcomes, of probability of an event as the proportion of
outcomes favourable to the event, and of the binomial distribution:
%
\begin{equation}
\label{bernoulli3} P(X=x) = \pmatrix{n
\cr
x} p^x(1-p)^{n-x},\qquad
x=0,1,2, \ldots, n,
\end{equation}
for the number $X$ of occurrences of the event in $n$ binomial trials.

Jacob Bernoulli's Theorem has two central features. The first is that
the greater the number of observations, the less the uncertainty. That
is: in a probabilistic sense later formalized as ``convergence in
probability'', relative frequencies of occurrence of an event in
independent repetitions of an experiment approach the probability of
occurrence of the event as sample size increases. It is in this guise
that Jacob Bernoulli's Theorem appears as the first \textit{limit theorem}
of probability theory, and in frequentist mathematical statistics as
the notion of a consistent estimator of a parameter (in this case the
parameter $p$).

The first central feature also reflects, as a mathematical theorem, the
empirically observed ``statistical regularity'' in nature, where
independent repetitions of a random experiment under uniform conditions
result in observed stability of large sample relative frequency of an event.

The second central feature, less stressed, is that Jacob Bernoulli's
Theorem is an \textit{exact result}. It is tantamount to obtaining a
sample size $n$ large enough for specified accuracy of approximation of
$p$ by the proportion $X/n$. The lower bound on such $n$
may depend on $p$, as it does in Jacob Bernoulli's Theorem, but even if
$p$ is known, the problem of determining the best possible lower bound
for $n$ for specified precision
is far from straightforward, as we shall demonstrate on one of
Bernoulli's examples.

Bernoulli's underlying motivation was, however, the approximation of an
\textit{unknown} $p$ by $X/n$ on the basis of repeated binomial sampling
(accumulation of evidence) to specific accuracy. We shall call this the
\textit{inversion} problem. It adds several layers of complexity to both features.

\subsection{Some background notes}
The present author's mathematical and historical interests have been
much in the direction of A.A. Markov, and Markov chains. It is
therefore a pleasure to have been asked to write this paper at the
Tricentenary for a journal which bears the name \textit{Bernoulli},
since A.A. Markov wrote an excellent summary of the history of the LLN
for the Bicentenary celebrations in St. Petersburg, Russia,
1913.\footnote{It is available in English as Appendix 1 of
\citet{ondar}.}

J.V. Uspensky's translation into Russian in 1913 of the fourth part of
\textit{Ars Conjectandi} (Bernoulli (\citeyear{bernoulli13},
\citeyear{bernoulli05})), where Jacob Bernoulli's Theorem occurs was
part of the St. Petersburg celebrations. Markov's paper and Uspensky's
translation are in \citet{bernoulli86}, a book prepared for the
First World Congress of the Bernoulli Society for Mathematical
Statistics and Probability held in Tashkent in 1986. It was one of the
present author's sources for this Tricentenary history, and includes a
long commentary by \citet{prokhorov}.

The Tricentenary of the death of Jacob Bernoulli was commemorated in
Paris in 2005 at a colloquium entitled \textit{L'art de conjecturer des
Bernoulli}. The proceedings have been published in the \textit{Journal
Electronique d'Histoire des Probabilit\'es et de la Statistique},
\textbf{2}, Nos. 1 and 1(b) (at
\href{http://www.jehps.net/}{www.jehps.net}). A number of celebrations
of the Tricentenary of publication of the \textit{Ars Conjectandi} are
scheduled for 2013, which is also the 250th anniversary of the first
public presentation of Thomas Bayes's work. This incorporates his
famous theorem, which plays an important role in our sequel. 2013 has
in fact been designated \textit{International Year of Statistics.}

Many of the sources below are available for viewing online, although
only few online addresses are specified in the sequel.

Titles of books and chapters are generally given in the original
language. In the case of Russian this is in English transliteration,
and with an English translation provided. English-language versions are
cited where possible. Quotations are in English translation. In French
only the first letter of titles is generally capitalized. German
capitalizes the first letters of nouns wherever nouns occur. Within
quotations we have generally stayed with whatever style the original
author had used.

\section{The Bernoullis and Montmort}

In 1687 Jacob Bernoulli (1654--1705) became Professor of Mathematics at
the University of Basel, and remained in this position until his death.

The title \textit{Ars Conjectandi} was an emulation of \textit{Ars
Cogitandi}, the title of the Latin version\footnote{Latin was then the
international language of scholarship. We have used ``Jacob'' as version
of the Latin ``Jacobus'' used by the author of \textit{Ars Conjectandi} for
this reason, instead of the German ``Jakob''.} of \textit{La Logique ou
l'Art de penser}, more commonly known as the Logic of Port Royal, whose
first edition was in 1662, the year of Pascal's death. Bienaym\'e
writes in 1843 (\citet{heydeseneta77}, p. 114) of Jacob Bernoulli:

\begin{quotation}
\noindent One reads on p. 225 of the fourth part of
his \textit{Ars
Conjectandi} that his ideas have been suggested to him, partially at
least, by Chapter 12 and the chapters following it, of \textit{l'Art de
penser}, whose author he calls \textit{magni acuminis et ingenii vir} [\textit{a
man of great acumen and ingenuity}] \ldots\ The final chapters contain
in fact elements of the calculus of probabilities, applied to history,
to medicine, to miracles, to literary criticism, to incidents in life,
etc., and are concluded by the argument of Pascal on eternal life.
\end{quotation}

The implication is that it was Pascal's writings which were the
influence. Jacob Bernoulli was steeped in Calvinism (although well
acquainted with Catholic theology). He was thus a firm believer in
predestination, as opposed to free will, and hence in determinism in
respect of ``random'' phenomena. This coloured his view on the origins
of statistical regularity in nature, and led to its mathematical formalization.

Jacob Bernoulli's \textit{Ars Conjectandi} remained unfinished in its
final (fourth) part, the \textit{Pars Quarta}, the part which contains the
theorem, at the time of his death. The unpublished version was reviewed
in the \textit{Journal des s\c{c}avans} in Paris, and the review
accelerated the (anonymous) publication of Montmort's \textit{Essay
d'analyse sur les jeux de hazard} in 1708.

Nicolaus Bernoulli (1687--1759)\footnote{For a substantial biographical
account, see \citet{csorgo}.} was a nephew to Jacob and Johann. His
doctorate in law at the University of Basel in 1709, entitled \textit{De
Usu Artis Conjectandi in Jure}, was clearly influenced by a direction
towards applications in the draft form of the \textit{Ars Conjectandi} of
Jacob. Nicolaus's uncle Johann, commenting in 1710 on problems in the
first edition of 1708 of Montmort, facilitated Nicolaus's contact with
Montmort, and seven letters from Nicolaus to Montmort appear in
\citet{montmort}, the second edition of the \textit{Essay}. The most
important of these as regards our present topic is dated Paris, 23
January, 1713. It focuses on a lower bound approximation to binomial
probabilities in the spirit of Jacob's in the \textit{Ars Conjectandi},
but, in contrast, for \textit{fixed} $n$.

As a specific illustration, Nicolaus asserts (\citet{montmort}, pp.
392--393), using modern notation, that if $X \sim B(14\,000, 18/35)$, then
%
\begin{equation}
\label{bernoulli4} 1-P(7037 \leq X \leq7363) \geq1/44.58 =0.0224316.
\end{equation}
\citet{laplace14}, p. 281, without mentioning Nicolaus anywhere, to
illustrate his own approximation, obtains the value $P(7037 \leq X \leq
7363) =0.994505$ so that $ 1 - P(7037 \leq X \leq7363) = 0.0056942$.
The statistical software package \textbf{R} which can calculate binomial
sum probabilities gives $P(X \leq7363) - P(X \leq7036)= 0.9943058$,
so that $ 1 - P(7037 \leq X \leq7363) = 0.0056942$.

The difference in philosophical approach is clear: finding $n$
sufficiently large for specific precision (Jacob), and finding the
degree of precision for given large $n$ (Nicolaus). While Nicolaus's
contribution is a direct bridge to the normal approximation of binomial
probabilities for large fixed $n$ as taken up by De Moivre, it does not
enhance the limit theorem direction of Jacob Bernoulli's Theorem, as De
Moivre's Theorem, from which the normal approximation for large $n$ to
binomial probabilities emerges, was to do.

After Paris, in early 1713 at Montmort's country estate, Nicolaus
helped Montmort prepare the second edition of his book
(\citet{montmort}), and returned to Basel in April, 1713, in time
to write a preface to \textit{Ars Conjectandi} which appeared in August
1713, a few months before Montmort's (\citeyear{montmort}).

Nicolaus, Pierre R\'emond de Montmort (1678--1719) and Abraham De
Moivre (1667--1754) were the three leading figures in what Hald
(\citeyear{hald07}) calls ``the great leap forward'' in stochastics,
which is how Hald describes the period from 1708 to the first edition
of De Moivre's (\citeyear{demoivre18}) \textit{Doctrine of Chances}.

In his preface to \textit{Ars Conjectandi} in 1713, Nicolaus says of the
fourth part that Jacob intended to apply what had been exposited in the
earlier parts of \textit{Ars Conjectandi}, to civic, moral and economic
questions, but due to prolonged illness and untimely onset of death,
Jacob left it incomplete. Describing himself as too young and
inexperienced to do this appropriately, Nicolaus decided to let the
\textit{Ars Conjectandi} be published in the form in which its author left
it. As \citet{csorgo} comments:

\begin{quotation} \noindent Jakob's programme, or dream rather, was
wholly impossible to accomplish in the eighteenth century. It is
impossible today, and will remain so, in the sense that it was
understood then.
\end{quotation}

\section{De Moivre}
De Moivre's motivation was to approximate sums of individual binomial
probabilities when $n$ is large, and the probability of success in a
single trial is $p$. Thus, when $ X \sim B(n,p) $. His initial focus
was on the symmetric case $p=1/2$ and large $n$, thus avoiding the
complication of approximating an asymmetric binomial distribution by a
symmetric distribution, the standard normal. In the English translation
of his 1733 paper (this is the culminating paper on this topic; a
facsimile of its opening pages is in \citet{stigler}, p.~74), De Moivre
(1738) praises the work of Jacob and Nicolaus Bernoulli on the summing
of several terms of the binomial $(a+b)^n $ when $n$ is large, which De
Moivre had already briefly described in his \textit{Miscellanea Analytica}
of 1730, but says:

\begin{quotation}\noindent \ldots\ yet some things were further
required; for what they have done is not so much an Approximation as
the determining of very wide limits, within which they demonstrated
that the sum of the terms was contained.
\end{quotation}

\noindent De Moivre's approach is in the spirit of Nicolaus
Bernoulli's, not least in that he seeks a result for large $n$, and
proceeds by approximation of individual terms.

As with Jacob Bernoulli's Theorem the limit theorem aspect of De
Moivre's result, effectively the Central Limit Theorem for the
standardized proportion of successes in $n$ binomial trials as $ n \to
\infty$, is masked, and the approximating value for sums of binomial
probabilities is paramount.

Nevertheless, De Moivre's results provide a strikingly simple, good,
and easy to apply approximation to binomial sums, in terms of an
integral of the normal density curve. He discovered this curve though
he did not attach special significance to it. \citet{stigler}, pp.
70--88 elegantly sketches the progress of De Moivre's
development. Citing \citet{stigler}, p. 81:

\begin{quotation}\noindent\ldots\ De Moivre had found an effective,
feasible way of summing the terms of the binomial.
\end{quotation}

\noindent It came about, in our view, due to two key components: what
is now known as Stirling's formula, and the practical calculation of
the normal integral.

\textit{De Moivre's} (\citeyear{demoivre33}) \textit{Theorem} may be stated as follows in modern
terms: the sum of the binomial terms
\[
\sum\pmatrix{n
\cr
x} p^xq^{n-x},
\]
where $0 <p =1-q <1$ over the range $|x-np| \leq s \sqrt{npq}$,
approaches as $n \to\infty$, for any $s>0$ the limit
%
\begin{equation}
\label{moivre1} \frac{1}{\sqrt{2 \pi}} \int_{-s}^{s}
e^{-z^2/2}\,dz.
\end{equation}

Jacob Bernoulli's Theorem as expressed by (\ref{bernoulli1}) \textit{follows as a Corollary.} This corollary is the LLN aspect of Jacob
Bernoulli's Theorem, and was the focus of De Moivre's application of
his result. It revolves conceptually, as does Jacob Bernoulli's
Theorem, around the mathematical formalization of statistical
regularity, which empirical phenomenon De Moivre attributes to:

\begin{quotation}\noindent\ldots\ that Order which naturally results
from ORIGINAL DESIGN.
\end{quotation}

\noindent(quoted by \citet{stigler}, p. 85). The theological
connotations of empirical statistical regularity in the context of free
will and its opposite, determinism, are elaborated in \citet{seneta03}.

De Moivre's (\citeyear{demoivre33}) result also gives an answer to estimating precision
of the relative frequency $X/n$ as an estimate of an \textit{unknown}
$p$, for given $n$; or of determining $n$ for given precision (the
inverse problem), in frequentist fashion, using the
inequality\footnote{See our Section 9.1.} \mbox{$p(1-p) \leq1/4$}.

De Moivre's results appeared in part in 1730 in his \textit{Miscellanea
Analytica de Seriebus et Quadraturis} and were completed in 1733 in his
\textit{Approximatio ad Summam Terminorum Binomii ${\overline{a+b}}|^n$ in
Seriem Expansi.}\footnote{${\overline{a+b}}|^n$ is De Moivre's
notation for $(a+b)^n$.} His \textit{Doctrine of Chances} of 1738 (2nd
ed.) contains his translation into English of the 1733 paper. There is
a short preamble on its p. 235, reproduced in \citet{stigler}, p. 74,
which states:

\begin{quotation}\noindent I shall here translate a Paper of mine which
was printed \textit{November} 12, 1733, and communicated to some Friends,
but never yet made public, reserving to myself the right of enlarging
my own thoughts, as occasion shall require.
\end{quotation}

In his \textit{Miscellanea Analytica}, Book V, De Moivre displays a
detailed study of the work of the Bernoullis in 1713, and distinguishes
clearly, on p. 28, between the approaches of Jacob in 1713 of finding
an $n$ sufficiently large for specified precision, and of Nicolaus of
assessing precision for fixed $n$ for the ``futurum probabilitate'',
thus alluding to the fact that the work was for a general, and to be
estimated, $p$.

The first edition of De Moivre's \textit{Doctrine of Chances} had appeared
in 1718. In this book there are a number of references to the work of
both Jacob and Nicolaus but only within a games of chance setting, in
particular to the work of Nicolaus as presented in \citet{montmort}, and
it is in this games of chance context that later French authors
generally cite the \textit{Doctrine of Chances}, characteristically giving
no year of publication.

The life and probabilistic work of De Moivre is throughly described
also in Schneider (\citeyear{schneider68}, \citeyear{schneider06}), and \citet{bellhouse}.

\section{Laplace. The inversion problem. Lacroix. The centenary}
In a paper of 1774 (Laplace (\citeyear{bernoulli86})) which \citet{stigler} regards as foundational for the
problem of predictive probability, Pierre Simon de Laplace (1749--1827)
sees that Bayes's
Theorem provides a means to solution of Jacob Bernoulli's \textit{inversion problem.} Laplace considers binomial trials with success
probability $x$ in each trial, assuming $x$ has uniform prior
distribution on $(0,1)$, and calculates the posterior distribution of
the success probability random variable $\Theta$ after observing $p$
successes and $q$ failures. Its density is:
%
\begin{equation}
\label{laplace1} \frac{\theta^p(1-\theta)^q}{\int_0^1\theta^p(1-\theta)^q \,d\theta} = \frac{(p+q+1)!}{p!q!}\theta^p(1-
\theta)^q
\end{equation}
and Laplace proves that for any given $w >0, \delta>0$
%
\begin{equation}
\label{laplace2} P\biggl(\biggl|\Theta- \frac{p}{p+q} \biggr| < w\biggr) > 1 - \delta
\end{equation}
for large $p, q$. This is a Bayesian analogue of Jacob Bernoulli's
Theorem, the beginning of Bayesian estimation of success probability of
binomial trials and of Bayesian-type limit theorems of LLN and Central
Limit kind. Early in the paper Laplace takes the \textit{mean}
%
\begin{equation}
\label{laplace3} \frac{p+1}{p+q+1}
\end{equation}
of the posterior distribution as his total [predictive] probability on
the basis of observing $p$ and~$q$, and (\ref{laplace3}) is what we now
call the Bayes estimator.

There is a brief mention in Laplace's paper of De Moivre's
\textit{Doctrine of Chances} (no doubt the 1718 edition) at the outset,
but in a context different from the Central Limit problem. There is no
mention of Jacob Bernoulli, Stirling (whose formula he uses, but which
he cites as sourced from the work of Euler), or Bayes. The paper of
1774 appears to be a work of still youthful exuberance.

In his preliminary \textit{Discours} to his \textit{Essai}, \citet{condorcet85},
p. viij, speaks of the relation between relative frequency and
probability, and has a footnote:

\begin{quotation}\noindent For these two demonstrations, \textit{see} the
third part of the \textit{Ars Conjectandi of Jacob Bernoulli}, a work full
of genius, and one of those of which one may regret that this great man
had begun his mathematical career so late, and whose death has too soon
interrupted.
\end{quotation}

\citet{lacroix}, who had been a pupil of J.A.N. de Caritat de Condorcet
(1743--1794), writes on p. 59 about Jacob Bernoulli's Theorem, and has
a footnote:

\begin{quotation} \noindent It is the object of the 4th Part of the
\textit{Ars Conjectandi}.
This posthumous work, published in 1713, already contains the principal
foundations of the philosophy of the probability calculus, but it
remained largely obscured until Condorcet recalled, perfected and
extended it.
\end{quotation}

Condorcet was indeed well-versed with the work of Jacob Bernoulli, and
specifically the \textit{Ars Conjectandi}, to which the numerous allusions
in the book of Bru--Crepel (\citet{condorcet94}) testify.

\citet{lacroix} may well be regarded as marking the \textit{first
Centenary} of Jacob Bernoulli's Theorem, because it gives a direct
proof and extensive discussion of that theorem. Subsequently, while the
name and statement of the theorem persist, it figures in essence as a
frequentist corollary to De Moivre's Theorem, or in its Bayesian
version, following the Bayesian (predictive) analogue of De Moivre's
Theorem, of Laplace (\citeyear{laplace14}, pp. 363 ff, Chapitre VI: \textit{De la
probabilit\'e des causes et des \'ev\'enemens futurs, tir\'ee des
\'ev\'enemens observ\'ees}), which is what the footnote of \citet{lacroix}, p. 295, cites at its very end.

The first edition of 1812 and the second edition of 1814 of Laplace's
\textit{Th\'eorie analytique des probabilit\'es} span the Centenary year
of 1813, but, as \citet{armatte} puts it, \citet{lacroix} served as an
exposition of the probabilistic
thinking of Condorcet and Laplace for people who would never go to the
original philosophical, let alone technical, sources of these
authors.\footnote{The influence of Lacroix's (\citeyear{lacroix}) book is particularly
evident in the subsequent more statistical direction of French
probability in the important book of \citet{cournot}, as the copious and
incisive notes of the editor, Bernard Bru, of its reprinting of 1984
make clear.}

Nevertheless, \citet{laplace14} is an outstanding epoch in the development
of probability theory. It connects well with what had gone before and
with our present history of the LLN, and mightily influenced the
future. Laplace's (\citeyear{laplace14}) p. 275 ff, Chapitre III, \textit{Des lois de
probabilit\'e, qui resultent de la multiplication ind\'efinie des
\'even\'emens} is frequentist in approach, contains De Moivre's
Theorem, and in fact adds a continuity correction term (p. 277):
%
\begin{equation}
\label{laplace4} P\bigl(|X - np| \leq t\sqrt{npq}\bigr) \approx\frac{1}{\sqrt{2\pi}} \int
_{-t}^{t} e^{-u^2/2} \,du + \frac{e^{-t^2/2}}{\sqrt{2 \pi npq}}.
\end{equation}
Laplace remarks that this is an approximation to $O(n^{-1})$ providing
$np$ is an integer,\footnote{See our Section 11.2 for a precise
statement.} and then applies it to Nicolaus Bernoulli's example (see
our {Section 2}). On p. 282 he inverts (\ref{laplace4}) to give an
interval for $p$ centred on ${\hat p}= X/n$, but the ends of the
interval still depend on the unknown $p$, which Laplace replaces by
${\hat p}$, since $n$ is large. This gives an interval of random
length, in fact a confidence interval in modern terminology, for $p$.

Neither De Moivre nor Stirling nor Nicolaus Bernoulli are mentioned
here. However in his \textit{Notice historique sur le Calcul des
Probabilit\'es}, pp. xcxix--civ, both Bernoullis, Montmort, De Moivre
and Stirling receive due credit. In particular a paragraph extending
over pages cij--ciij refers to a later edition (1838 or 1856,
unspecified) of De Moivre's \textit{Doctrine of Chances} specifically in
the context of De Moivre's Theorem, in both its contexts, that is (1)
as facilitating a proof of Jacob Bernoulli's Theorem; and (2) as:

\begin{quotation} \noindent\ldots\ an elegant and simple expression
that the difference between these two ratios will be contained within
the given limits.
\end{quotation}

Finally, of relevance to our present theme is a subsection (pp. 67--70)
entitled: \textit{Th\'eor\`emes sur le developpement en s\'eries des
fonctions de plusieurs variables.} Here Laplace considers, using their
generating functions, sums of independent integer-valued but not
necessarily identically distributed random variables, and obtains a
Central Limit Theorem. The idea of inhomogeneous sums and averages
leads directly into subsequent French (Poisson) and Russian (Chebyshev)
directions.

\section{Poisson's Law of Large Numbers and Chebyshev}
\subsection{Poisson's Law}
The major work in probability of Sim\'eon Denis Poisson (1781--1840)
was his book\footnote{The digitized version which I have examined has
the label on the cover: ``From the Library of J.V. Uspensky, Professor
of Mathematics at Stanford, 1929--1947'' and is from the Stanford
University Libraries. Uspensky plays a major role in our account.} of
1837: \textit{Recherches sur la probabilit\'e}. It is largely a treatise
in the tradition of, and a sequel to, that of his great predecessor
Laplace's (\citeyear{laplace14}) \textit{Th\'eorie analytique} in its emphasis on the
large sample behaviour of averages. The theorem of Jacques Bernouilli
[sic] [Jacob Bernoulli] is mentioned in 5 places, scattered over pp.
iij to p. 205. Laplace receives multiple mentions on 16 pages. Bayes,
as ``Blayes'', is mentioned twice on just one page, and in connection
with Laplace. What follows is clearly in the sense of Laplace, with the
prior probability values for probability of success in binomial trials
being determined by occurrence of one of a range of ``causes''.
Condorcet is mentioned twice, and Pascal at the beginning, but there is
no mention of Montmort, let alone Nicolaus Bernoulli, nor of De Moivre.

The term \textit{Loi des grands nombres} [\textit{Law of Large Numbers}] appears
for the first time in the history of probability on p. 7 of \citet{poisson}, within the statement;

\begin{quotation} \noindent Things of every kind of nature are subject
to a universal law which one may well call \textit{the Law of Large
Numbers}. It consists in that if one observes large numbers of events
of the same nature depending on causes which are constant and causes
which vary irregularly$, \ldots,$ one finds that the proportions of
occurrence are almost constant \ldots
\end{quotation}

There are two versions of a LLN in Poisson's treatise. The one most
emphasized by him has at any one of $n$ binomial trials, each of a
fixed number $a$ of causes operate equiprobably, that is, with
probability $1/a$, occurrence of the $i$th cause resulting in observed
success with probability $p_i, i =1,2, \ldots, a$. Thus in each of
$n$ independent trials the probability of success is $ {\bar p}(a) =
\sum_{i=1}^n p_i/a$. So if $X$ is the number of successes in $n$
trials, for sufficiently large $n$,
\[
P\biggl(\biggl|\frac{X}{n} - {\bar p}(a) \biggr| > \epsilon\biggr) < Q
\]
for any prespecified $\epsilon, Q$. \citet{poisson} proved this
directly, not realizing that it follows directly from Jacob Bernoulli's Theorem.

The LLN which \citet{poisson} considered first, and is now called
Poisson's Law of Large Numbers, has probability of success in the $i$th
trial fixed, at $p_i, i = 1, 2, \ldots, n$. He showed that
\[
P\biggl(\biggl|\frac{X}{n} - {\bar p}(n) \biggr| > \epsilon\biggr) < Q
\]
for sufficiently large $n$, using Laplace's Central Limit Theorem for
sums of non-identically distributed random variables. The special case
where $p_i = p, i=1,2, \ldots$ gives Jacob Bernoulli's Theorem, so
Poisson's LLN is a genuine generalization.

Inasmuch as $ {\bar p}(n)$ itself need not even converge as $ n \to
\infty$, Poisson's LLN displays as a primary aspect \textit{loss of
variability} of proportions $X/n$ as $ n \to\infty$, rather than a
tendency to \textit{stability}, which Jacob Bernoulli's Theorem
established under the restriction $p_i = p$.

\subsection{Chebyshev's thesis and paper}

The magisterial thesis, \citet{chebyshev45}, at Moscow University of
Pafnutiy Lvovich Chebyshev (1821--1894), begun in 1841 and defended in
1846, but apparently published in Russian in 1845, was entitled \textit{An
Essay in Elementary Analysis of the Theory of Probabilities.}\footnote
{I have consulted a reprinting in \citet{chebyshev55}, pp. 111--189.} It
gives as its motivation, dated 17 October (o.s) 1844 (\citet{chebyshev55},
pp. 112--113):

\begin{quotation} \noindent To show without using transcendental
analysis the fundamental theorems of the calculus of probabilities and
their main applications, to serve as a support for all branches of
knowledge, based on observations and evidence \ldots
\end{quotation}

Dominant driving forces for the application of probability theory in
Europe, Great Britain, and the Russian Empire in those times were
\textit{retirement funds and insurance}\footnote{\citet{laplace14} had
devoted an extensive part of the applications to investigations of life
tables and the sex ratio, and in France De Moivre's work was largely
known for his writings on annuities.} and Russian institutions such as
the Yaroslavl Demidov Lyc\'ee, within the Moscow Educational Region,
had no textbooks. Such a textbook was to involve only ``elementary
methods''.

As a consequence, Chebyshev's magisterial dissertation used no
calculus, only algebra, with what would have been integrals being sums
throughout, but was nevertheless almost entirely theoretical, giving a
rigorous analytical discussion of the then probability theory, with a
few examples. Throughout, the quantity $e^{-x^2}$ figures prominently.
The dissertation concludes with a table of what are in effect tail
probabilities of the standard normal distribution. Much of the thesis
is in fact devoted to producing these very accurate tables (correct to
7 decimal places) by summation.

Laplace's (\citeyear{laplace14}) Chapitre VI, on predictive probability, is adapted by
Chebyshev to the circumstances.
In Laplace's writings, the prior distribution is envisaged as coming
about as the result of ``causes'', resulting in corresponding values
being attached to the possible values in $(0,1)$ of a success
probability, the attached value depending on which ``cause'' occurs. If
causes are deemed to be ``equiprobable'', the distribution of success
probability is uniform in $(0,1)$.

Chebyshev stays ``discrete'', so, for example, he takes $\frac{i}{s}, i
=1,2, \ldots, s-1$ as the possible values (the sample space) of the
prior probability in $(0,1)$ of $ (s-1)$ equiprobable causes, the
probability of each of the causes being $\frac{1}{s-1}$. Thus if $r$
occurrences of an event $E$ (``success'') are observed in $n$ trials,
the posterior distribution is given by:
%
\begin{equation}
\label{causes3} \frac{ {n \choose r}(i/s)^r(1-(i/s)^{n-r}}{\sum_{i=1}^{s-1}{n \choose
r} (i/s)^r(1-(i/s))^{n-r}}.
\end{equation}

Examples are also motivated by \citet{laplace14}, who in the same chapter
begins Section~28, p. 377, with the following:

\begin{quotation}\noindent It is principally to births that the
preceding analysis is applicable.
\end{quotation}

Chebyshev's (\citeyear{chebyshev45}) thesis concludes Section 26,which is within Chapter
IV, with:

\begin{quotation} \noindent Investigations have shown that of 215 599
newborns in France 110 312 were boys.\footnote{I could not find this
data in \citet{laplace14}, although more extensive data of this kind is
treated there.}
\end{quotation}

He then calculates that the probability that the posterior random
variable $\Theta$ satisfies
\[
P( 0.50715 \leq\Theta\leq0.51615) = 0.99996980
\]
by taking $r/n =110312/215599 = 0.511653579$, and using (in modern notation):
\[
\Theta\sim{\mathcal{N}} { \biggl(}\frac{r}{n}, \frac{(r/n)(1 -
(r/n))}{n}{
\biggr) }
\]
and his tables of the standardized normal random variable. (Using the
statistical software
\textbf{R} for the standard normal variable gives 0.9999709.)

Chebyshev is clearly well acquainted with not only the work of Laplace,
but also the work of De Moivre, Bayes and Stirling, although he cites
none of these authors explicitly. Jacob Bernoulli's Theorem is
mentioned at the end of Chebyshev's (\citeyear{chebyshev45}) thesis,
Section~20, where he proceeds to obtain as an approximation to the
binomial probability:
\[
P_{\mu, m} = \frac{\mu!}{m! (\mu- m)!} p^m(1-p)^{\mu- m}
\]
the expression
\[
\frac{1}{\sqrt{2 \pi p(1-p)\mu}}e^{- {\frac{z^2}{2}} /2p(1-p)\mu}
\]
using the
(Stirling) approximation
$ x! = \sqrt{2\pi} x^{x+ 1/2} e^{-x}$ which he says is the ``form
usually used in probability theory''. But he actually obtains bounds for
$x!$ directly.\footnote{See our Section 9.2.} Much of Chapter III is in
fact dedicated to finding such bounds, as he says at the outset to this chapter.

Such careful bounding arguments (rather than approximative asymptotic
expressions) are characteristic of Chebyshev's work, and of the Russian
probabilistic tradition which came after him. This is very much in the
spirit of the bounds in Jacob Bernoulli's Theorem.

Poisson's (\citeyear{poisson}) \textit{Recherches sur la probabilit\'e} came to
Chebyshev's attention after the publication of \citet{chebyshev45}, but
the subtitle of \citet{chebyshev46} suggests that the content of \citet{chebyshev46}, motivated by Poisson's LLN was used in the defense of the
dissertation in 1846 (\citet{bernstein45}).

The only explicit citation in \citet{chebyshev46} is to \citet{poisson},
Chapitre IV, although Jacob Bernoulli's Theorem is acknowledged as a
special case of Poisson's LLN. In his Section 1 Chebyshev says of
Poisson's LLN:

\begin{quotation} \noindent All the same, no matter how ingenious the
method utilized by the splendid geometer, it does not provide bounds on
the error in this approximative analysis, and, in consequence of this
lack of degree of error, the derivation lacks appropriate rigour.
\end{quotation}

\citet{chebyshev46} in effect repeats his bounds for the homogeneous case
($p_i = p= 1,2, \ldots, n $) of binomial trials which he dealt with in
\citet{chebyshev45}, Section 21, to deal with the present inhomogeneous
case. He also uses generating functions for sums in the manner of
\citet{laplace14}.

Here is his final result, where as usual $X$ stands for the number of
successes in $n$ trials, $p_i$ is the probability of success in the
$i$th trial, and $p = \frac{\sum_{i=1}^{n} p_i}{n}$.
%
\begin{equation}
\label{cheb1} P\biggl(\biggl|\frac{X}{n} - p\biggr| \geq z\biggr) \leq Q \mbox{ if } n
\geq\max{ \biggl\{} \biggl( \frac{\log[Q {\frac{z}{1-p}} \sqrt{ \frac
{1-p -z}{p+z}}]}{\log H} \biggr), \biggl( \frac{\log[Q {\frac
{z}{p}} \sqrt{ \frac{p -z}{1-p+z}}]}{\log H_1}
\biggr) { \biggr\}}
\end{equation}
where:
%
\begin{equation}
\label{cheb2} H = \biggl( \frac{p}{p+z} \biggr)^{p+z} \biggl(
\frac{1-p}{1-p-z} \biggr)^{1-p-z},\qquad H_1 = \biggl(
\frac{p}{p-z} \biggr)^{p-z} \biggl(\frac
{1-p}{1-p+z}
\biggr)^{1-p+z}.
\end{equation}

Structurally, (\ref{cheb1}), (\ref{cheb2}) are very similar to Jacob
Bernoulli's expressions in his Theorem, so it is relevant to compare
what they give in his numerical example when $z = 1/50, p =30/50 =0.6,
Q = 1/1001= 0.000999001$. The answer appears to be $ n \geq12241.293$,
i.e., $ n \geq12242$.

In spite of the eminence of the journal (Crelle's) in which \citet{chebyshev46} published, and the French language in which he wrote, the paper
passed unnoticed among the French mathematicians, to whom what we now
call Poisson's LLN remained an object of controversy. In his
historically important follow-up to the Laplacian analytical tradition
of probability (see \citet{bru12}), \citet{laurent} gives a
proof of Poisson's Law of Large Numbers. He uses characteristic
functions, and gives a careful consideration of the error, and hence of
convergence rate. However \citet{sleshinsky}, in his historical
Foreward, claims Laurent's proof contains an error which substantially
alters the conclusion on convergence rate. \citet{laurent} cites a
number of Bienaym\'e's papers, but does not appear to use the simple
proof of Poisson's LLN which follows from the Bienaym\'e--Chebyshev
Inequality, which by 1873 had been known for some time.

\section{Bienaym\'e and Chebyshev}

\subsection{Bienaym\'e's motivation}

The major early work of Iren\'ee Jules Bienaym\'e (1796--1878): \textit{De
la dur\'ee de la vie en France} (1837), on the accuracy of life tables
as used for insurance calculations, forced the abandonment of the
Duvillard table in France in favour of the Deparcieux table. He was
influenced in the writing of this paper not least by the demographic
content of Laplace's \textit{Th\'eorie analytique}.

Bienaym\'e, well aware of \citet{poisson} vehemently disapproved of the
term ``Law of Large Numbers'' (\citet{heydeseneta77}, Section 3.3),
thinking that it did not exist as a separate entity from Jacob
Bernoulli's Theorem, not understanding the version of Poisson's Law
where a \textit{fixed} probability of success, $p_i$ is associated with
the $i$-th trial, $ i =1,2, \ldots\,$.
As a consequence of his misunderstanding, in 1839 (\citet{heydeseneta77}, Section 3.2) Bienaym\'e proposes a scheme of variation of
probabilities (that is, of ``genuine'' inhomogeneity of trials, as
opposed to the other version of Poisson's Law which does not differ
from Jacob Bernoulli's) through a principle of \textit{dur\'ee des
causes}
[\textit{persistence of causes}]. Suppose there are $a$ causes, say $c_1, c_2,
\ldots, c_a$, the $i$-th cause giving rise to probability $p_i, i = 1,
2, \ldots, a$ of success. Each cause may occur equiprobably for each
one of $m$ sets of $n$ trials; but once chosen it \textit{persists} for
the whole set. The case $n=1$ is of course the ``other'' Poisson scheme
which is tantamount to Jacob Bernoulli's sampling scheme with success
probability ${\bar p}(a)$.

For his scheme of $N=mn$ trials Bienaym\'e writes down a Central Limit
result with correction term to the normal integral in the manner of
Laplace's version of De Moivre's Theorem for Bernoulli trials, to which
Bienaym\'e's result reduces when $n=1$.

Schemes of $m$ sets of $n$ binomial trials underlie Dispersion Theory,
the study of homogeneity and stability of repeated trials, which was a
predecessor of the ``continental direction of statistics''. Work on
Dispersion Theory proceeded through Lexis, Bortkiewicz and Chuprov; and
eventually, through the correspondence between Markov and Chuprov,
manifested itself in another branch of the evolutionary tree of the LLN
of repeated trials founded on Jacob Bernoulli's Theorem. (See \citet{heydeseneta77}, Chapter 3.)

\subsection{The Bienaym\'e--Chebyshev Inequality}

Bienaym{\'e} (\citeyear{biename53}) shows mathematically that for the sample mean ${\bar
X}$ of independently and identically distributed random variables whose
population mean is $\mu$ and variance is $\sigma^2$, so $ E {\bar X} =
\mu, \operatorname{Var} {\bar X} = \sigma^2/n$, then for any $t>0$:
%
\begin{equation}
\label{inequality1} \operatorname{Pr}\bigl(({\bar X} - \mu)^2 \geq t^2
\sigma^2\bigr) \leq1/\bigl(t^2n\bigr).
\end{equation}
The proof which Bienaym\'e uses is the simple one we use in the
classroom today to prove the inequality by proving that for any
$\epsilon>0$, providing $E X^2 < \infty$, and $\mu= EX$:
%
\begin{equation}
\label{inequality2a} \operatorname{Pr}\bigl(|X - \mu|
\geq\epsilon\bigr) \leq(\operatorname{Var} X)/\epsilon^2.
\end{equation}
This is commonly referred to in probability theory as the Chebyshev
Inequality, and less commonly as the Bienaym\'e--Chebyshev Inequality.
If the $X_i, i = 1,2, \ldots$ are independently but not necessarily
identically distributed, and $ S_n = X_1+X_2+ \cdots+X_n$, putting $X
=S_n$ in (\ref{inequality2a}), and using the Bienaym\'e equality $ \operatorname{Var}
S_n = \Sigma_{i=1}^n \operatorname{Var}X_i, $ (\ref{inequality2a}) reads:
%
\begin{equation}
\label{inequality3} \operatorname{Pr}\bigl(|S_n - ES_n| \geq\epsilon\bigr) \leq
\Biggl(\Sigma_{i=1}^n \operatorname{Var} X_i\Biggr)\Big/
\epsilon^2.
\end{equation}
This inequality was obtained by \citet{chebyshev67} for discrete random
variables and published simultaneously in French and Russian.
\citet{biename53} was reprinted immediately preceding the French
version in Liouville's journal. In 1874 Chebyshev wrote:

\begin{quotation} \noindent The simple and rigorous demonstration of
Bernoulli's law to be found in my note entitled: \textit{Des valeurs
moyennes}, is only one of the results easily deduced from the method of
M. Bienaym\'e, which led him, himself, to demonstrate a theorem on
probabilities, from which Bernoulli's law follows immediately \ldots
\end{quotation}

\noindent Actually, not only the limit theorem aspect of Jacob
Bernoulli's Theorem is covered by the Bienaym\'e--Chebyshev Inequality,
but also the inversion aspect\footnote{Using $p(1-p) \leq1/4$.} even
for unspecified~$p$.

Further, \citet{chebyshev74} formulates as: ``the method of
Bienaym\'e'' what later became known as the method of moments.
Chebyshev (\citeyear{Che87}) used this method to prove the first
version of the Central Limit Theorem for sums of independently but not
identically distributed summands; and it was quickly taken up and
generalized by Markov. Markov and Liapunov were Chebyshev's most
illustrious students, and Markov was ever a champion of Bienaym\'e as
regards priority of discovery. See \citet{heydeseneta77}, Section
5.10, for details, and
\citet{seneta84}\footnote{\href{http://www.maths.usyd.edu.au/u/eseneta/TMS_9_37-77.pdf}{http://www.maths.usyd.edu.au/u/eseneta/TMS\_9\_37-77.pdf}.}
for a history of the
Central Limit problem in pre-Revolutionary Russia.

\section{Life tables, insurance, and probability in Britain. De Morgan}

From the mid 1700s, there had been a close association between games of
chance and demographic and official statistics with respect to
calculation of survival probabilities from life tables. Indeed games of
chance and demographic statistics were carriers of the nascent
discipline of probability. There was a need for, and activity towards,
a reliable science of risk based on birth statistics and life tables by
insurance companies and superannuation funds
(\citet{heydeseneta77}, Sections 2.2--2.3). De Moivre's
(\citeyear{demoivre25}) \textit{Annuities upon Lives} was a
foremost source in England.

John William Lubbock (1803--1865) is sometimes described as ``the
foremost among English mathematicians in adopting Laplace's doctrine of
probability''. With John Elliott Drinkwater (Later Drinkwater-Bethune)
(1801--1859), he published anonymously a 64 page elementary treatise on
probability (\citet{lubbockdrw}). Lubbock's slightly
younger colleague, Augustus De Morgan (1806--1871), was making a name
for himself as mathematician, actuary and academic. The paper of
\citet{lubbock} attempts to address and correct current shortcomings of
life tables used at the time in England. He praises Laplace's
\textit{Th\'eorie analytique} in respect of its Bayesian approach, and
applies this approach to multinomial distributions of observations, to
obtain in particular probabilities of intervals symmetric about the
mean via a normal limiting distribution.

Lubbock very likely used the 1820 edition of the \textit{Th\'eorie
analytique}, since his colleague De Morgan was in 1837 to review this
edition.
De Morgan's
chief work on probability was, consequently, a lengthy article in the
\textit{Encyclopedia Metropolitana} usually cited as of 1845, but
published as separatum in 1837 (\citet{demorgan37}).
This was primarily a summary, simplification and clarification of many
of Laplace's derivations. It was the first full-length exposition of
Laplacian theory and the first major work in English on probability theory.

An early footnote (p. 410) expresses De Morgan's satisfaction that in
Lubbock and Drinkwater-Bethune (\citeyear{lubbockdrw}) there is a collection in English,
``and in so accessible a form'' on ``problems on gambling which usually
fill works on our subject'', so he has no compunction in throwing most
of these aside ``to make room for principles'' in the manner of, though
not necessarily in the methodology of, Laplace. There is no mention of
De Moivre or Bayes.

On p. 413, Section 48, which is on ``the probability of future events
from those which are past'', De Morgan addresses the same problem as
\citet{lubbock}. Using the multinomial distribution for the prior
probabilities, he calculates the posterior distribution by Bayes's
Theorem, and this is then used to find the joint distribution from
further drawings. Stirling's formula (with no attribution) is
introduced in Section 70. A discussion of the normal approximation to
the binomial follows in Section~74, pp. 431--434. We could find no
mention as such of Jacob Bernoulli's Theorem or De Moivre's Theorem.
Section~74 is concluded by Nicolaus Bernoulli's example, which is taken
directly from Laplace: with success probability $18/35$, and $14\,000$
trials, $P(7200-163 \leq X \leq7200+163)$ is considered, for which De
Morgan obtains 0.99433 (a~little closer to the true value 0.99431 than
Laplace). Section 77, p. 434, addresses ``the inverse question'' of
prediction given observations and prior distribution.

De Morgan (\citeyear{demorgan38}) published \textit{An Essay on Probabilities}, designed
for the use of actuaries. The book, clearly written and much less
technical than \citet{demorgan37} remained widely used in the insurance
industry for many years. It gave an interesting perception of the
history up to that time, especially of the English contributions. On
pp. v--viii De Morgan says:

\begin{quotation} \noindent At the end of the seventeenth century, the
theory of probability was contained in a few isolated problems, which
had been solved by Pascal, Huyghens, James Bernoulli, and others. \ldots
\ Montmort, James Bernoulli, and perhaps others, had made some slight
attempts to overcome the mathematical difficulty; but De Moivre, one of
the most profound analysts of his day, was the first who made decided
progress in the removal of the necessity for tedious operations \ldots\
when we look at the intricate analysis by which Laplace obtained the
same [results], \ldots\
De Moivre nevertheless did not discover the inverse method. This was
first used by the Rev. T. Bayes, \ldots\ Laplace, armed with the
mathematical aid given by De Moivre, Stirling, Euler and others, and
being in possession of the inverse principle already mentioned,
succeeded \ldots\ in \ldots\ reducing the difficulties of calculation\
\ldots\ within the reach of an ordinary arithmetician \ldots\ for the
solution of all questions in the theory of chances which would
otherwise require large numbers of operations. The instrument employed
is a table (marked Table I in the Appendix to this work), upon the
construction of which the ultimate solution of every problem may be
made to depend.
\end{quotation}

Table I is basically a table of the normal distribution.

\section{The British and French streams continue}
\subsection{Boole and Todhunter}
George Boole (1815--1864), a prot\'eg\'e of De Morgan, in his book
\citet{boole}, introduced (p.~307) what became known as Boole's
Inequality, which was later instrumental in coping with statistical
dependence in Cantelli's (\citeyear{cantelli17})\vadjust{\goodbreak} pioneering treatment
of the Strong Law of Large Numbers (\citet{seneta92}). Boole's
book contains one of the first careful treatments of hypothesis testing
on the foundation of Bayes's Theorem.\footnote{According to
\citet{todhunter} there is no difference in essence between the
1814 2nd and the 1820 3rd editions.} (\citet{riceseneta}). Boole
does not appear to pay attention to Jacob Bernoulli's Theorem, nor does
he follow Laplacian methods although he shows respect for De Morgan, as
one who:

\begin{quote}\noindent has most fully entered into the spirit of
Laplace.
\end{quote}

He objects to the uniform prior to express ignorance, and to inverse
probability (Bayesian) methods in general, particularly in regard to
his discussion of Laplace's Law of Succession. Boole is kinder to
\citet{poisson}, whom he quotes at length at the outset of his Chapter
XVI: \textit{On the Theory of Probabilities}. His approach to this theory
is in essence set-theoretic, in the spirit of formal logic.

The history of the theory of probability upto and including Laplace,
and even some later related materials, appears in the remarkable book
of \citet{todhunter} which is still in use today. A whole chapter
entitled: Chapter VII. James Bernoulli (Sections 92--134, pp. 56--77)
addresses the whole of the \textit{Ars Conjectandi}. Sections 123--124,
devoted to Jacob Bernoulli's Theorem, begin with:

\begin{quote} \noindent The most remarkable subject contained in the
fourth part of the \textit{Ars Conjectandi} is the enunciation of what we
now call \textit{Bernoulli's Theorem}.
\end{quote}

The theorem is enunciated just as Bernoulli described it; of how large
$N (n)$ is to be to give the specified precision. Section 123 ends with:

\begin{quote} \noindent James Bernoulli's demonstration of this result
is long but perfectly satisfactory \ldots We shall see that James
Bernoulli's demonstration is now superseded by the use of Stirling's
Theorem.
\end{quote}

In Section 124, Todhunter uses Jacob Bernoulli's own examples,
including the one we have cited (``for the odds to be 1000 to 1'').
Section 125 is on the \textit{inversion} problem: given the number of
successes in $n$ trials, to determine the precision of the estimate of
the probability of success. Todhunter concludes by saying that the
inversion has been done in two ways,

\begin{quote} \noindent by an inversion of James Bernoulli's Theorem,
or by the aid of another theorem called Bayes's theorem; the results
approximately agree. See Laplace \textit{Th\'eorie Analytique \ldots}
pages 282 and 366.
\end{quote}

Section 135 concludes with:

\begin{quotation}\noindent The problems in the first three parts of the
\textit{Ars Conjectandi} cannot be considered equal in importance or
difficulty to those which we find investigated by Montmort and De
Moivre; but the memorable theorem in the fourth part, which justly
bears its author's name, will ensure him a permanent place in the
history of the Theory of Probability.
\end{quotation}

Important here is Todhunter's view that Jacob Bernoulli's \textit{proof}
has been superseded. Only the \textit{limit theorem aspect} is being
perceived, and that as a corollary to De Moivre's Theorem, although in
this connection and at this point De Moivre gets no credit, despite
Laplace's (\citeyear{laplace14}) full recognition for his theorem.

\subsection{Crofton and Cook Wilson}
Todhunter's limited perception of Jacob Bernoulli's Theorem as only a
limit theorem, with Stirling's Theorem as instrument of proof in the
manner of De Moivre--Laplace, but without mention of De Moivre, became
the standard one in subsequent British probability theory. In his
\textit{Encyclopaedia Britannica} article in the famous 9th edition,
\citet{crofton} constructs such a proof (pp. 772--773), using his
characteristically geometric approach, to emphasize the approximative
use of the normal integral to calculate probabilities, and then
concludes with:

\begin{quotation} \noindent Hence \textit{it is always possible to
increase the number of trials till it becomes certainty that the
proportion of occurrences of the event will differ from $p$} (\textit{its
probability on a single trial}) \textit{by a quantity less than any
assignable.} This is the celebrated theorem given by James Bernoulli in
the \textit{Ars Conjectandi.} (See Todhunter's \textit{History}, p. 71.)
\end{quotation}

Then Crofton presents the whole issue of Laplace's predictive approach
as a consequence of Bayes's Theorem in Section 17 of \citet{crofton}
(pp. 774--775), using a characteristically
geometric argument, together with Stirling's Theorem.

Crofton's general francophilia is everywhere evident; he had spent a
time in France. His concluding paragraph on p. 778, on literature,
mentions De Morgan's \textit{Encyclopaedia Metropolitana} presentation,
Boole's book with some disparagement, and a number of French language
sources, but he refers:

\begin{quotation}\noindent\ldots\ the reader, \ldots\ above all, to
the great work of Laplace, of which it is sufficient to say that it is
worthy of the genius of its author -- the \textit{Th\'eorie analytique des
probabilit\'es, \ldots}
\end{quotation}

There is a certain duality between De Moivre, a Protestant refugee from
Catholic France to Protestant England, and Crofton, an Anglo-Irish
convert to Roman Catholicism in the footsteps of John Henry (Cardinal)
Newman, himself an author on probability, and an influence on Crofton,
as is evident from its early paragraphs, in \citet{crofton}.

Crofton's (\citeyear{crofton}) article was likely brought to the
attention \cite{seneta12} of John Cook Wilson (1849--1915), who in
\citet{cookwilson} developed his own relatively simple proof of
the limit aspect of ``James Bernoulli's Theorem''. He uses domination
by a geometric progression. His motivation is the simplification of
Laplace's proof as presented in Todhunter (\citeyear{todhunter},
Section~993). There is no mention of De Moivre, and dealings with the
normal integral are avoided. An interesting feature is that Cook Wilson
considers \textit{asymmetric} bounds for the deviation $\frac{X}{n} -
p$, but he does eventually resort to limiting arguments using
Stirling's approximation, so the possibility of an \textit{exact}
bounding result in the Bernoulli and Bienaym\'e--Chebyshev style is
lost.

\subsection{Bertrand}
In the book of \citet{bertrand} in the French stream, Chapitre IV
contains a proof of De Moivre's Theorem, and mentions both De Moivre
and Stirling's Theorem, but there seems to be no mention of ``Jacques
Bernoulli'' in the chapter content, nor a statement of his Theorem, let
alone a proof. Chapitre V of \citet{bertrand} has two ``demonstrations''
which Bertrand (\citeyear{bertrand}, p. 101) describes only in its limit aspect.
Bertrand first shows that if $ X_i, i =1, \ldots, n $ are independently
and identically distributed, and $ EX_1^2 < \infty$, then $ \operatorname{Var} {\bar
X}= \frac{\operatorname{Var} X_1}{n} \to0, n \to\infty$ and then simply applies this
to the case when $P(X_1=1) = p, P(X_1=0) = q =1-p$. There is no mention
of the Bienaym\'e--Chebyshev Inequality or its authors. That $ \operatorname{Var}
{\bar X} \to0 $ is deemed sufficient for ``convergence'' one might
charitably equate to foreshadowing convergence in mean square. In the
final section of Chapitre V, Section 80, p. 101, \citet{bertrand}
asserts that he will give a demonstration to the theorem of Bernoulli
even simpler than the preceding. What follows is a demonstration that
for $ \{0, 1 \} $ random variables $X_i, i = 1, \ldots, n, E| {\bar X}
- p | \to0, n \to\infty$ without the need to calculate $ E|
\Sigma_{i=1}^n X_i - np|. $ The reader will see that this actually
follows easily from $\operatorname{Var} {\bar X} \to0$. The ``exact aspect'' of Jacob
Bernoulli's theorem has disappeared.

\subsection{K. Pearson}
In a perceptive paper written in that author's somewhat abrasive style,
Karl \citet{pearson} restores credit to De Moivre for his achievement,
and refocuses (p. 202) on the need, surely appropriate for a
mathematical statistician, to obtain a better expression for sample
size, $n$, needed for specified accuracy.

In his Section 2, Pearson reproduces the main features of Jacob
Bernoulli's proof, and shows how the normal approximation to the
binomial in the manner of De Moivre can be used to determine $n$ for
specified precision \textit{if $p$ is known}. In Section 3, Pearson
tightens up Bernoulli's proof keeping the same expressions for $p =\frac
{r}{r+s}$ and $\epsilon=\frac{1}{r+s}, $ by using\vspace*{1pt} a geometric series
bounding procedure and then Stirling's Theorem. There is no mention of
Cook Wilson's (\citeyear{cookwilson}) work. Recall that if $p \neq\frac{1}{2}$ one
problem with the normal approximation to the normal is that asymmetry
about its mean of the binomial is not reflected in the normal. Thus in
considering
%
\begin{equation}
\label{pearson1} \frac{c}{c+1} < P\biggl( \biggl| \frac{X}{n} - p\biggr| \leq
\epsilon\biggr) = P(X \leq np +n \epsilon) -P (X < np - n \epsilon)
\end{equation}
involves binomial tails of differing probability size.

A commensurate aspect in \citet{pearson} is the treatment of the tails
of the binomial distribution individually. The approximation is
remarkably good, giving for Bernoulli's example where $ r=30, s=20, p =
\frac{3}{5}, c =1000, \epsilon= \frac{1}{50}$ the result $ n_0(\epsilon
, c) \geq6502$, which is almost the same as for the normal
approximation to the binomial (6498). The reason is similar: the use
of the De Moivre--Stirling approximation for $x!$, and the fact that
$p=0.6$ is close to $p=0.5$, which is the case of symmetric binomial
(when it is known that a correction for continuity such as Laplace's
with the normal probability function gives very accurate results).
Pearson does not attempt the inversion (that is, the determination of
$n$ when $p$ is not known) in Jacob Bernoulli's example.

\section{Sample size and emerging bounds}
\subsection{Sample size in Bernoulli's example}
For this classical example when $p=0.6$, referring to (\ref{pearson1}),
we seek the smallest $n$ to satisfy
%
\begin{equation}
\label{pearson2} 0.9990009999= \frac{1000}{1001} < P(X \leq0.62n ) - P (X <
0.58n)
\end{equation}
where $X \sim B(n,0.6)$. Using \textbf{R}, $n =6491$ on the right hand
side gives $0.9990126$, while $n= 6490$ gives $0.9989679$, so the
minimal $n$ which will do is $6491$, providing the algorithm in \textbf{R}
is satisfactory.

Chebyshev's (\citeyear{chebyshev46}) inequality for inhomogeneous binomial trials, when
applied to a homogeneous situation, gives, as we have seen, a much
sharper ``exact'' result for minimal $n$ (namely, $n \geq12242$) for
$p=0.6$ than Bernoulli's, but, like Jacob Bernoulli's, was incapable of
explicit algebraic inversion when $p$ was unknown.

In his monograph (in the 3rd edition, \citet{markov13}, this is on p. 74)
Markov uses the normal approximation with known $p=0.6$ in Bernoulli's
example to obtain that $ n \geq6498 $ is required for the specified
accuracy.\footnote{Actually, Markov uses $0.999$ in place of
$0.9990009999$.}

In the tradition of Chebyshev, and in the context of his controversy
with Nekrasov (see our Section 10.1), \citet{markov99} had developed a
method using continued fractions to obtain tight \textit{bounds} for
binomial probabilities when $p$ is known and $n$ is also prespecified.
The method is described and illustrated in \citet{uspensky}, pp. 52--56.
On p. 74, \citet{markov13} argues that the upper bound $0.999$ on accuracy
is likely to hold also for $n$ not much greater than the approximative
$6498$ which he had just obtained, say $n=6520$. On pp. 161--165 he
verifies this, showing that the probability when $n=6520$ is between
$0.999028$ and $0.999044$. Using \textbf{R} the true value is 0.9990309.
Thus Markov's procedure is an exact procedure for inversion when $p$
and accuracy are prespecified, once one has an approximative lower
bound for $n$. One could then proceed experimentally, as we have done
using \textbf{R}, looking for the smallest $n$.

To effect ``approximative'' inversion \textit{if we did not know the value
of $p$}, to get the specified accuracy of the \textit{estimate} of $p$
presuming $n$ would still be large, we could use De Moivre's Theorem
and the ``worst case'' bound $p(1-p) \leq\frac{1}{4}$, to obtain
\[
n \geq\frac{z_0^2}{4 \epsilon^2} = 0.25(3.290527)^2 (50)^2
=6767.23 \geq6767
\]
where $P(|Z| \leq z_0) = 0.999001. $ Again the result $6767$ is good
since $p=0.6$ is not far from the worst case value $p=0.5$. The now
commonly used substitution of the estimate $\hat p$ from a preliminary
performance of the binomial experiment in place of $p$ in $p(1-p)$ (and
this is implicit in Laplace's use of his add-on correction to the
normal to effect inversion) would improve the inversion result.

\subsection{Improving Stirling's approximation and the normal approximation}
The De Moivre--Laplace approximative methods are based on Stirling's
approximation for the factorial. They can be refined by obtaining \textit{bounds} for the factorial. Such bounds were already present in an
extended footnote in \citet{chebyshev46}:
%
\begin{equation}
T_0 x^{x+\frac{1}{2}}e^{-x} < x! < T_0
x^{x+\frac{1}{2}}e^{-x+ \frac{1}{12x}}
\end{equation}
where $T_0$ is a positive constant. This was later refined\footnote{For
a history see \citet{boudin}, pp. 244--251: Note II. Formule de
Stirling.} to
%
\begin{equation}
x! =\sqrt{2 \pi}x^{x+\frac{1}{2}}e^{-x+ \frac{1}{12x+\theta}}
\end{equation}
where $0<\theta<1$.

It was therefore to be expected that De Moivre's Theorem could be made
more precise by producing bounds. \citet{delavallee1} in the
second\vspace*{1pt} of two papers (see Seneta (\citeyear{seneta01a})) considers the sum $ P= \Sigma_{x}
{n \choose x } p^xq^{n-x} $ over the range $ |x- (n+1)p +\frac
{1}{2}| < (n+1) {l }$ for arbitrary fixed ${l}$, and obtains the bounds
for $P$ in terms of the normal integral. A bound for minimal sample
size $n$ required for specified accuracy of approximation could be
determined, at least when $p$ was known. Although this work seems to
have passed largely unnoticed, it presages the return of ``exact
methods'' via bounds on the deviation of normal approximation to the
binomial. These bounds imply a \textit{convergence rate} of $O(n^{-1/2})$.

A cycle of related bounding procedures is initiated in \citet{bernstein11}. He begins by saying that he has not found a rigorous estimate
of the accuracy of the normal approximation (``Laplace's formula'') to $
P(|X - np|< z \sqrt{np(1-p)}). $ He illustrates his own investigations
by showing when: $p=1/2$, $n$ is an odd number, and $ 1/2 + z\sqrt {(n+1)/2}$ is an integer, that:
%
\begin{equation}
\label{bernstein1} P\biggl(\biggl|X-\frac{n}{2}\biggr| \leq z\sqrt{\frac{n+1}{2}
}\biggr) > 2\Phi(z \sqrt{2}) -1
\end{equation}
where $\Phi(z) = P(Z \leq z)$ is the cumulative distribution function
of a standard normal variable $ Z \sim\mathcal{N}(0,1)$.
He illustrates in the case $z=2.25, n =199$, so that the right-hand
side\footnote{Using \textbf{R}, the left-hand side is $0.9989406$.} of
(\ref{bernstein1}) is $0.9985373$. This value is thus a lower bound for
$P(77 \leq X \leq122)$. \citet{bernstein11} then inverts, by finding that
if $\Phi(z_0) = 0.9985373$, then this value normal approximation
corresponds to $P(77.05 \leq X \leq121.9) = P(78 \leq X \leq121)$.
This testifies not only to the accuracy of ``Laplace's formula'' as
approximation, but also to the sharpness of Bernstein's bound even for
moderate size of $n$, albeit in the very specific situation of $p=
1/2$.

The accuracy of Laplace's formula became a central theme in Bernstein's
subsequent probabilistic work. There is a strong thematic connection
between Bernstein's striking work on \textit{probabilistic methods} in
approximation theory in those early pre-war years to about 1914, in
Kharkov, and De La Vall\'ee Poussin's approximation theory. See our
Section 11.3.

\section{The Russian stream. Statistical dependence and Bicentenary
celebrations}

\subsection{Nekrasov and Markov}

\citet{nekrasov98a} is a summary paper containing no proofs. It is
dedicated to the memory of Chebyshev, on account of Nekrasov's
continuation of Chebyshev's work on Central Limit theory in it. The
author, P.A. Nekrasov (1853--1924), attempted to use what we now call
the method of saddle points, of Laplacian peaks, and of the Lagrange
inversion formula, to establish, for sums of independent
non-identically distributed lattice random variables, what
are now standard local and global limit theorems of Central Limit
theory for large deviations. A follow-up paper, \citet{nekrasov98b} dealt
exclusively with binomial trials.

Markov's (\citeyear{markov98}) first rigorization, within correspondence with A.V.
Vasiliev (1853--1929), of Chebyshev's version of the Central Limit
Theorem, appeared in the Kazan-based journal edited by Vasiliev. The
three papers of 1898 mark the beginning of two bitter controversies
between Nekrasov and Markov, details of whose technical and personal
interaction are described in \citet{seneta84}.

Nekrasov's writings from about 1898 had become less mathematically
focused, partly due to administrative load, and partly due to his use
of statistics as a propagandist tool of state and religious authority
(Tsarist government and the Russian Orthodox Church).

In a long footnote (\citet{nekrasov}, pp. 29--31), states ``Chebyshev's
Theorem'' as follows:
If $X_1, X_2, \ldots,X_n$ are independently distributed and $ {\bar
X}_n = (X_1 + X_2+ \cdots+ X_n)/n $ then
\[
P\bigl(|{\bar X}_n - E{\bar X}_n | < \tau
\sqrt{g_n}\bigr) \geq1 - \frac{1}{n
\tau^2},
\]
where $\tau$ is a given positive number, and
\[
g_n =\frac{ \Sigma_{i=1}^n \operatorname{Var} X_i}{n}.
\]
He adds that if $ \tau(=\tau_n) $ can be chosen so that $ \tau_n \sqrt {g_n} \to0 $ while simultaneously $ n \tau^2_n \to\infty$, then $
{\bar X}_n - E{\bar X}_n $ converges to $0$. This comment encompasses
the LLN in its general form at the time.

Nekrasov says (\citet{seneta84}) that he has examined the ``theoretical
underpinnings of Chebyshev's Theorem'', and has come to the [correct]
conclusion that
if in the above $g_n$ is defined as $ g_n = n \operatorname{Var} {\bar X}_n$, the
inequality continues to hold. Now, in general,
\[
\operatorname{Var} {\bar X}_n = \frac{ \Sigma_{i=1}^n \operatorname{Var} X_i + 2 \Sigma_{i < j}
\operatorname{Cov}(X_i, X_j)}{n^2},
\]
so the \textit{original expression} for $g_n$ results from just ``pairwise
independence'' (i.e. all $ \operatorname{Cov}(X_i, X_j) =0, {i < j})$.

Hence if under merely \textit{pairwise
independence} $\tau_n$ can be chosen so that $ \tau_n \sqrt{g_n} \to0
$, the LLN will hold. Thus (under this latter condition) pairwise
independence is \textit{sufficient} for the LLN. This is an important
advance: a first step for the LLN to hold under a condition weaker than
the hitherto assumed mutual independence.

\citet{nekrasov}, p. 29, then boldly states that Chebyshev's
Theorem\footnote{Among Russian authors of the time, such as Nekrasov
and A.A. Chuprov, the LLN itself is often called Chebyshev's Theorem,
with no distinction of the Bienaym\'e--Chebyshev Inequality
from its application.} attains its ``\textit{full force}'' under the
condition that the $ X_i, i \geq1$ are \textit{pairwise independent}.

There is little doubt from his context that Nekrasov is asserting that
pairwise independence is \textit{necessary} for the LLN to hold. Markov
saw that this was not correct, and proceeded to construct a
counterexample: the first ``Markov'' chain (\citet{seneta96}).

All that was needed was an example of dependence where $g_n = n \operatorname{Var}
{\bar X}_n $ and $\tau_n$ is such that $\tau_n^2 g_n \to0$, while
simultaneously $n \tau_n^2 \to0$.
Still publishing in his friend Vasiliev's journal, \citet{markov06},
Sections 2 and 5, does almost precisely this (\citet{seneta96}, Section
5). He constructs in fact, as his \textit{general scheme of dependent
variables} $ \{ X_n \},\break n \geq1$, a finite Markov chain, which he
takes to be time homogeneous with all transition probabilities $
p_{ij}$ strictly positive.
He shows that $ E X_n $ has a limit, $a$, as $ n \to\infty$ (that
limit is in fact the mean of the limiting-stationary distribution), and
then uses the Bienaym\'e--Chebyshev Inequality to show that $ P(| {\bar
X}_n - a| \geq\epsilon) \to0$ as $ n \to\infty$.

The last sentence of \citet{markov06} reads (without mention of Nekrasov):

\begin{quotation} \noindent Thus, independence of quantities does not
constitute a necessary condition for the existence of the law of large
numbers.
\end{quotation}

\subsection{Chuprov and Markov}

In his book (Chuprov (\citeyear{chuprov10})) which was a pioneering and fundamental influence on
statistics in the Russian Empire, \textit{Essays on the Theory of
Statistics}, A.A. Chuprov (1874--1926) speaks of the following
ideological conflict in thinking about the LLN, especially in Russia.
The LLN is a mathematical theorem. It \textit{reflects} an empirical fact:
observed long-term stability of the proportion of successes in
independent binomial trials, \textit{but does not explain its cause}
(\citet{seneta03}).

The consequent Markov--Chuprov Correspondence (\citet{ondar}) lasted from
2 November 1910 to about 27 February 1917, and marks the coming
together of probability theory and statistics (in the form of
Dispersion Theory) into mathematical statistics in the Russian Empire.
The Correspondence itself was largely concerned, in the tradition of
the Lexis--Bortkiewicz theory, with the study of variants of the
empirical dispersion coefficient.
A full account of Dispersion Theory more broadly is in \citet{heydeseneta77}, Section 3.4.

The Correspondence refers frequently to the work on Dispersion Theory
of Bortkiewicz, a~Russian expatriate of Polish ethnicity, in Germany.
Bortkiewicz's (\citeyear{bort}) book, \textit{Das Gesetz der kleinen Zahlen}. [\textit{The
Law of Small Numbers.}] (LSN) has earned itself a niche in the history
of mathematical statistics. The name is clearly in contrast to
Poisson's LLN, but what precisely it describes is unclear. A relatively
recent study is by \citet{quineseneta}.

The time-span of the Markov--Chuprov Correspondence encompassed \textit{the
Bicentenary} of Jacob Bernoulli's Theorem,
and this did not pass unnoticed. In \citet{ondar}, p. 65, Letter No. 54
(a letter from Markov to Chuprov, 15 January, 1913):

\begin{quotation} \noindent Firstly, do you know: the year 1913 is the
two hundredth anniversary of the law of large numbers (\textit{Ars
Conjectandi, 1713}), and don't you think that this anniversary should
be commemorated in some way or other? Personally I propose to put out a
new edition of my book, substantially expanded. But I am raising the
question about a general celebration with the participation of as large
a number of people and institutions as possible.
\end{quotation}

Then in \citet{ondar}, p. 69, Letter No. 56 (a letter from Markov to
Chuprov, 31 January, 1913) anticipates a view which
\citet{pearson} was to express even more forcefully later:

\begin{quotation} \noindent\ldots\ Besides you and me, it was
proposed to bring in Professor A.V. Vasiliev \ldots\  Then it was
proposed to translate only the fourth chapter of \textit{Ars Conjectandi};
the translation will be done by the mathematician Ya.V. Uspensky, who
knows the Latin language well, and it should appear in 1913. Finally, I
propose to do a French translation of the supplementary articles in my
book with a short foreword about the law of large numbers, as a
publication of the Academy of Sciences. All of this should be scheduled
for 1913 and a portrait of J. Bernoulli will be attached to all the
publications.

In connection with your idea about attracting foreign scholars, I
cannot fail to note that the French mathematicians, following the
example of Bertrand, do not wish to know what the theorem of Jacob
Bernoulli is. In Bertrand's \textit{Calcul des Probabilit\'es} the fourth
chapter is entitled ``Th\'eor\`eme de Jacques Bernoulli''. the fifth is
entitled ``D\'emonstration \'el\'ementaire du th\'eor\`eme de Jacques
Bernoulli''. But neither a strict formulation of the theorem nor a
strict proof is given. \ldots\ Too slight a respect for the theorem of
J. Bernoulli is also observed among many Germans. It has reached a
point that a certain
Charlier \ldots\ showed a complete lack of familiarity with the
theorem.
\end{quotation}

A.V. Vasiliev (1853--1929) was both mathematician and social activist, a
Joe Gani of his time and place. Among his interests was the history of
mathematics in Russia. He was at Kazan University over the years
1874--1906, as Professor apart from the early years; and at St.
Petersburg--Petrograd\footnote{The name of the city was changed from
Sankt Peterburg (sometimes written Sanktpeterburg) to Petrograd during
World War I, then to Leningrad after the Bolshevik Revolution, and is
now Sankt Peterburg again.} University 1907--1923. Only slightly older
than Markov (1856--1922), he was instrumental in fostering Markov's work
over the period 1898--1906 in ``his'' Kazan journal: the \textit{Izvestiia}
of the Physico-Mathematical Society of Kazan University.

\subsection{The Bicentenary in St. Petersburg}

The talks were in the order: Vasiliev, Markov, Chuprov. Markov's talk,
published in Odessa in 1914, is in \citet{ondar} in English translation,
as Appendix 3, pp. 158--163. Chuprov's much longer talk, also published
in 1914, is Appendix 4, pp. 164--181.

Vasiliev (\citet{vasiliev}) presents a brief summary of the whole
proceedings, and then in his Sections I, II, III, the respective
contributions of the three speakers, in French (in a now-electronically
accessible journal). He describes the respective topics as: Vasiliev:
Some questions of the theory of probabilities upto the theorem of
Bernoulli; Markov: The Law of Large Numbers considered as a collection
of mathematical theorems; Chuprov: The Law of Large Numbers in
contemporary science.

The content of his own talk is described as giving a historical
perception of the development of two fundamental notions: mathematical
probability, that is \textit{a priori}; and empirical probability, that is
\textit{a posteriori}. Vasiliev summarizes Markov's talk well, especially
the early part which contrasts Jacob Bernoulli's \textit{exact} results
with the \textit{approximative procedures} of De Moivre and Laplace, which
use the limit normal integral structure to determine probabilities.
Markov mentions Laplace's second degree correction, and also comments
on the \textit{proof} of Jacob Bernoulli's Theorem in its limit aspect by
way of the De Moivre--Laplace ``second limit theorem''.

Markov goes on to discuss Poisson's LLN as an approximative procedure
``\ldots\ not bounding the error in an appropriate way'', and
continues with Chebyshev's (\citeyear{chebyshev46}) proof in Crelle's journal. He then
summarizes the Bienaym\'e--Chebyshev interaction in regard to the
Inequality and its application; and the evolution of the method of
moments. He strains to remain scrupulously fair to Bienaym\'e, while
according Chebyshev credit, in a way adopted subsequently in
Russian-language historiography, for example by S.N. \citet{bernstein45}.
Markov concludes as follows, in a story which has become familiar.

\begin{quotation} \noindent\ldots\ I return to Jacob Bernoulli. His
biographers recall that, following the example of Archimedes he
requested that on his tombstone the logarithmic spiral be inscribed
with the epitaph ``Eadem mutato resurgo''. \ldots\ It also expresses
Bernoulli's hope for resurrection and eternal life. \ldots\ More than
two hundred years have passed since Bernoulli's death but he lives and
will live in his theorem.
\end{quotation}

Chuprov's bicentennial talk contrasts two methods of knowledge: the
study of the individual entity, and the study of the collective via
averages, ``based on the Law of Large Numbers''. As an illustration of
the success of the latter, Mendel's laws of heredity are cited. What
seems to be the essence here is the goodness of fit to a probability
model of repeated statistical observations under uniform conditions
(statistical sampling).

Markov perceives in such arguments the vexed question of statistical
regularity being interpreted as the LLN, which to him (and to us) is a
mathematical theorem which, under specific mathematical conditions,
only reflects statistical regularity, and does not explain it. He
writes, obviously miffed, to Chuprov immediately after the meeting
(\citet{ondar}, Letter No. 62, 3~December 1913):

\begin{quotation} \noindent In your talk statistics stood first and
foremost, and applications of the law of large numbers were advanced
that seem questionable to me. By subscribing to them I can only weaken
that which for me is particularly dear: the rigor of judgements I
permit. \ldots\ Your talk harmonized beautifully with A.V. Vasiliev's
talk but in no way with mine. \ldots\ I had to give my talk since the
200th anniversary of a mathematical theorem was being celebrated, but I
do not intend to publish it and I do not wish to.
\end{quotation}

Chuprov's paper as a whole is scholarly and interesting, for example,
also mentioning Brown (of Brownian motion), and the Law of Small
Numbers as related to the Poisson distribution, and to Abb\'e and Bortkiewicz.

As anticipated in Markov's letter (\citet{ondar}, No. 56, 31 January
1913) the translation from Latin into Russian by J.V. Uspensky was
published in 1913, edited, and with a Foreword, by Markov, and with the
now-usual portrait of Jacob Bernoulli. They are reproduced in
\citet{bernoulli86}, as is Markov's talk at the Bicentenary meeting. Additionally,
to celebrate the Bicentenary Markov published in 1913 the 3rd
substantially expanded edition of his celebrated monograph \textit{Ischislenie Veroiatnostei} [\textit{Calculus of Probabilities}], complete with
the portrait of, Jacob Bernoulli. The title page is headed:

\begin{quotation}
\noindent{\textbf{K 200 lietnemu iubileiu zakona bol'shikh chisel.}} [To
the 200th-year jubilee of the law of large numbers.]
\end{quotation}
with the title \textit{Ischislenie Veroiatnostei} below, with other information.
For the portrait of Jacob Benoulli following the title page, Markov, at
the conclusion of his Preface, expresses his gratitude to the chief
librarian of Basel University, Dr. Carl Christoph Bernoulli.

\section{\texorpdfstring{Markov (\protect\citeyear{markov13}), Markov's Theorems, Bernstein and Uspensky}{Markov (1913), Markov's Theorems, Bernstein and
Uspensky}}
\subsection{\texorpdfstring{Markov (\protect\citeyear{markov13}) and Markov's Theorems}{Markov (1913) and Markov's Theorems}}

In this 3rd Bicentenary edition, \citet{markov13}, Chapter III (pp.
51--112), is titled \textit{The Law of Large Numbers}, and Chapter IV (pp.
113--171) is titled \textit{Examples of various methods of calculation of
probabilities.} Chapters are further subdivided into numbered but
untitled subsections. The innovation which Markov feels most important,
according to his Preface, are the appendices (pp. 301--374):

\begin{quotation}\noindent Application of the method of mathematical
expectations -- the method of moments -- to the proof of the second
limit theorem of the calculus of probabilities.
\end{quotation}
The first appendix is titled \textit{Chebyshev's inequalities and the
fundamental theorem} and the second \textit{Theorem on the limit of
probability in the formulations of Academician A.M. Liapunov.} In a
footnote to the latter, Markov finally attributes the proof of De
Moivre's Theorem to \citet{demoivre30}.\footnote{Hitherto his
attributions had been to Laplace.}

What is of specific interest to us is what has come to be known as
Markov's Inequality\footnote{\citet{bernstein27} (1934), p. 101 and p. 92,
respectively, and then \citet{uspensky}, p. 182, call it Chebyshev's
(resp. Tshebysheff's) Lemma.}: for a non-negative random variable $U$
and positive number $u$:
%
\begin{equation}
\label{markovineq} P(U \geq u) \leq\frac{E(U)}{u}
\end{equation}
which occurs as a Lemma on pp. 61--63. It is then used to prove (\ref
{inequality2a}), which Markov calls the Bienaym\'e--Chebyshev
Inequality, on pp. 63--65, in what has become the standard modern
manner, inherent already in Bienaym\'e's (\citeyear{biename53}) proof.

Section 16 (of Chapter III) is entitled \textit{The Possibility of
Further Extensions.} It begins on p. 75 and on p. 76 Markov asserts
that
%
\begin{equation}
\label{markov1} \frac{\operatorname{Var}( S_n)}{n^2} \to0 \mbox{ as } n \to\infty
\end{equation}
is sufficient for the WLLN to hold,
for arbitrary summands $\{X_1, X_2, \ldots\}$. Thus the assumption of
independence is dropped, although the assumption of finite individual
variances is still retained. In the Russian literature, for example in
\citet{bernstein27}, p. 177, this is called \textit{Markov's Theorem.} We
shall call it Markov's Theorem 1.

Amongst the innovations in this 3rd edition which Markov actually
specifies in his Preface is an advanced version of the WLLN which came
to be known also as \textit{Markov's Theorem}, and which we shall call
Markov's Theorem 2. We state it in modern form:
%
\begin{equation}
\label{wlln1} \frac{S_n}{n} - E { \biggl(}\frac{S_n}{n}{ \biggr)}
\stackrel{p} {\to} 0
\end{equation}
where $S_n =\sum_{i=1}^n X_i$ and the $ \{X_i, i=1,2, \ldots\}$ are
independent and satisfy $E(|X_i|^{1+\delta}) < C < \infty$ for some
constants $\delta>0$ and $C$. The case $\delta= 1$ came to be known
in Russian-language literature as \textit{Chebyshev's Theorem.} Markov's
Theorem 2 thus dispenses with the need for finite variance of summands
$X_i$, but retains their independence. It occurs in the same Section 16
of Chapter III of \citet{markov13}, specifically on pp. 83--88.

Markov's publications of 1914 strongly reflect his apparent background
reading activity in preparation for the Bicentenary. In particular, a
paper entitled \textit{O zadache Yakova Bernoulli} [\textit{On the problem of Jacob
Bernoulli}] can be found in \citet{markov51}, pp. 509--521. In this paper
in place of what Markov calls the approximate formula of De Moivre:
\[
\frac{1} {\sqrt{\pi}}\int_z^{\infty}
e^{-z^2}\,dz \mbox{ for } P(X > np + z \sqrt{2npq})
\]
he derives the expression
\[
\frac{1} {\sqrt{\pi}}\int_z^{\infty}
e^{-z^2}\,dz + \frac
{(1-2z^2)(p-q)e^{-z^2}}{6\sqrt{2npq \pi}}
\]
which Markov calls Chebyshev's formula. This paper of Markov's clearly
motivated \citet{uspensky} to ultimately resolve the issues through the
component (\ref{uspensky3}) of Uspensky's expression.

\subsection{Bernstein and Uspensky on the WLLN}

Markov died in 1922 well after the Bolshevik seizure of power, and it
was through the 4th (posthumous) edition of \textit{Ischislenie
Veroiatnestei} (\citet{markov24}) that his results were publicized and
extended, in the first instance in the Soviet Union due to the
monograph S.N. \citet{bernstein27}. The third part of this book (pp.
142--199) is titled \textit{The Law of Large Numbers} and consists of three
chapters: Chapter 1: \textit{Chebyshev's inequality and its consequences.}
Chapter 2: \textit{Refinement of Chebyshev's Inequaliity.} and Chapter 3:
\textit{Extension of the Law of Large Numbers to dependent quantitities.}
Chapter 3 begins with Markov's Theorem 1, and continues with study of
the effect of specific forms of correlation between the summands
forming $S_n$. In Chapter 1, on p. 155
Bernstein mentions Markov's Theorem 2 as a result of the ``deceased
Academician A.A. Markov'' and adds ``The reader will find the proof in
the textbook of A.A. Markov''. A proof is included in the second
edition, Bernstein (\citeyear{bernstein34}), in which the three chapters in the third
part are almost unchanged from \citet{bernstein27}.

\citet{bernstein24} returned to the problem of accuracy of the normal
approximation to the binomial via bounds. He showed that there exists
an $\alpha (|\alpha| \leq1)$ such that $ P= \Sigma_{x} {n \choose x
} p^xq^{n-x} $ summed over $x$ satisfying
%
\begin{equation}
\biggl|x- np -\frac{t^2}{6}(q-p)\biggr| < t \sqrt {npq} + \alpha \mbox{ is } {
\frac{1}{\sqrt{2\pi}}\int_{-t}^{t} e^{-u^2/2}\,du}
+ 2\theta e^{-(2npq)^{1/3}}
\end{equation}
where $|\theta| <1$ for any $n, t$, provided $t^2/16 \leq npq \geq
365$. The tool used, perhaps for the first time ever, was what came to
be known as \textit{Bernstein's Inequality}:
%
\begin{equation}
\label{bernstein2} P(V>v) \leq\frac{E(e^{V\epsilon})}{e^{v\epsilon}}
\mbox{ for any $\epsilon>0$,
which follows from } P(U>u) \leq\frac{E(U)}{u}
\end{equation}
namely Markov's Inequality (called Chebyshev's Lemma by Bernstein). It
holds for any random variable $V$, on substituting $ U= e^{V\epsilon},
u = e^{v\epsilon}$. If $E(e^{V\epsilon})< \infty$ the bound is
particularly effective for a non-negative random variable $V$ such as
the binomial, since the bound may be tightened by manipulating
$\epsilon$. In connection with a discussion of (\ref{bernstein2}),
\citet{bernstein27}, pp. 231--232 points out that, consequently the
ordinary (uncorrected) normal integral approximation thus gives
adequate accuracy when $npq$ is of size several hundred, but in cases
where great accuracy is not required, $npq \geq30$ will do. However,
our interest in (\ref{bernstein2}) is in its nature as an \textit{exact}
result and in the suggested \textit{rate of convergence}, $O(n^{-1})$, to
the limit in the WLLN which the bounds provide.

The entire issue was resolved into an ultimate exact form, under the
partial influence of the extensive treatment of the WLLN in Bernstein's
(\citeyear{bernstein27}) textbook, by Uspensky (\citeyear{uspensky}, Chapter VII, p. 130) who showed that
$P$ taken over the usual range $ t_1 \sqrt{npq} \leq x-np \leq t_2
\sqrt{npq}$ for any real numbers $t_1 < t_2$, can be expressed
(provided $npq \geq25)$ as:
%
\begin{eqnarray}
\label{uspensky1}
{\frac{1}{\sqrt{2\pi}}\int_{t_1}^{t_2}
e^{-u^2/2}\,du} &+& \frac{(1/2 -
\theta_1)e^{-t_1^2/2} +(1/2 - \theta_2)e^{-t_2^2/2}}{\sqrt{2\pi npq}}
\\
\label{uspensky3} &+&\frac{(q-p)\{ (1- t_2^2)e^{-t_2^2/2}-(1-
t_1^2)e^{-t_1^2/2}\}}{6\sqrt{2\pi npq}} + \Omega,
\end{eqnarray}
where
$\theta_2 = np +t_2\sqrt{npq} - [np +t_2\sqrt{npq}], \theta_1 = np
-t_1\sqrt{npq} - [np -t_1\sqrt{npq}]$, and
\[
|\Omega| < \frac{0.20 + 0.25|p-q|}{npq} + e^{-3\sqrt{npq}/2}.
\]
The symmetric case then follows by putting $t_2 = - t_1 =t$,
so the ``Chebyshev'' term vanishes. When both $np$ and $ t \sqrt{npq} $
are integers, $\theta_1 = \theta_2 =0$, reducing the correction term in
(\ref{uspensky1}) to Laplace's $ e^{-t^2/2}/\sqrt{2\pi npq}$. But in
any case, bounds which are within O$(n^{-1})$ of the true value are
thus available.

Uspensky's (\citeyear{uspensky}) book carried Markov's theory to the English-speaking
countries. \citet{uspensky} cites \citet{markov24} and \citet{bernstein27} in
his two-chapter discussion of the LLN. Markov's Theorem 2 is stated and
proved in Chapter X, Section 8. Presumably the second (\citeyear{bernstein34}) edition of
Bernstein's textbook was not available to Uspensky due to circumstances
mentioned below. On the other hand in \citet{uspensky} the ideas in the
proof of Markov's Theorem 2 are used to prove the now famous
``Khinchin's Theorem'', an ultimate form of the WLLN. For independent
identically distributed (i.i.d.) summands, Khinchin (\citet{khinchin29})
showed that the existence of a finite mean, $ \mu= E X_i$, is
sufficient for (\ref{wlln1}). Finally, \citet{uspensky}, pp. 101--103,
proves the Strong Law of Large Numbers (SLLN) for the setting of
Bernoulli's Theorem, and calls this strengthening ``Cantelli's
Theorem'', citing one of the two foundation papers (\citet{cantelli17}) in
the history of the SLLN.

On the other hand, \citet{bernstein34}, in his third part has an
additional Chapter 4: \textit{Statistical probabilities, average values
and the coefficient of dispersion.} It begins with a Bayesian \textit{inversion} of Jacob Bernoulli's Theorem, proved under a certain
condition on the prior (unconditional) distribution of the number of
``successes'', X, in $n$ trials. The methodology uses Markov's
Inequality applied to $P  {(}(\Theta- \frac{X}{n})^4 > w^4  {|}
\Theta {)} $ and in the classical case of a uniform prior
distribution over $(0,1)$ of the success probability $\Theta$ gives for
any given $w >0$
%
\begin{equation}
\label{laplacebernstein2} P {\biggl(} \biggl|\Theta- \frac{X}{n} \biggr| < w {\Big|} X= m {
\biggr)}> 1 - \frac
{3(n_0+1)}{16nw^4n_0},
\end{equation}
for $n > n_0$
and $ m =0,1,\ldots, n$. This should be compared with (\ref{laplace2}).

\citet{bernstein34} also has 4 new appendices. The 4th of these (pp.
406--409) is titled: \textit{A Theorem Inverse to Laplace's Theorem.} This
is the Bayesian inverse of De Moivre's Theorem, with an arbitrary prior
density, and convergence to the standard normal integral as $ m, n \to
\infty$ providing $\frac{m}{n}$ behaves appropriately. A version of
this theorem is now called the Bernstein-von Mises Theorem, although
this attribution is not quite appropriate. After Laplace the
multivariate extension of Laplace's inversion is actually due to his
disciplele Bienaym\'e in 1834, and is called by von Mises in 1919 the
``Second Fundamental Theorem'' (the first being the CLT). Details are
given in Section 5.2 of \citet{heydeseneta77}.

The books of both Bernstein and Uspensky are very much devoted to
Markov's work, and Bernstein's also emphasizes and publicizes Markov
chains. Several sections of Bernstein's textbook in its 1946 4th
edition, such as the 4th appendix, are included in Bernstein's (\citeyear{bernstein64})
collected works and have not been published separately.

\subsection{Biographical notes on Bernstein and Uspensky}
Bernstein and Uspensky played parallel and influential roles in
publicizing and extending Markov's work, especially on the LLN.
These roles were conditioned by their background. To help understand,
we sketch these backgrounds. Uspensky's story is almost unknown.

Sergei Natanovich Bernstein (1880--1968) was born in Odessa in the
then Russian Empire. Although his father was a doctor and university
lecturer, the family had difficulties since it was Jewish. On
completing high school, Bernstein went to Paris for his mathematical
ediucation, and defended a doctoral dissertation in 1904 at the
Sorbonne. He returned in 1905 and taught at Kharkov University from
1908 to 1933. In the spirit of his French training and following a
Chebyshevian theme, in the years preceding the outbreak of World War I
he followed \citet{bernstein11} by a number of articles on approximation
theory. These included the famous paper of 1912 which presented a
probabilistic proof of Weierstrass's Theorem, and introduced what we
now call Bernstein Polynomials. A~prize-winning paper which contains
forms of inverse theorems and Bernstein's Inequality, arose out of a
question posed by De La Vall\'ee Poussin.\footnote{See History of
Approximation Theory (HAT) at
\url{http://www.math.technion.ac.il/hat/papers.html}.}

After the Bolshevik Revolution during 1919--1934 Kharkov (Kharkiv in
Ukrainian) was the capital of the Ukrainian SRS. Bernstein became
Professor at Kharkov University and was active in the Soviet
reorganization of tertiary institutions as National Commissar for
Education, when the All-Ukrainian Scientific Research Institute of
Mathematical Sciences was set up in 1928. In 1933 he was forced to move
to Leningrad, where he worked at the Mathematical Institute of the
Academy of Sciences. He and his wife were evacuated to Kazakhstan
before Leningrad was blockaded by German armies from September, 1941 to
January, 1943. From 1943 he worked at the Mathematical Institute in
Moscow. Further detail may be found in Seneta (\citeyear{seneta01b}).

\citet{bernstein64} is the 4th volume of the four volume collection of his
mathematical papers. His continuing interest in the accuracy of the
normal distribution as approximation to the binomial probabilities
developed into a reexamination in a new light of the main theorems of
probability, such as their extension to dependent summands. The idea of
martingale differences
appears in his work, which is perhaps best known for his extensions of
the Central Limit Theorem to ``weakly dependent random variables''. His
was a continuing voice of reason in the face of Stalinist interference
in mathematical and biological science. A fifth edition of his textbook
never appeared. It was stopped when almost in press because of
prevailing ideology.

J.V. Uspensky, the translator of the 4th part of \textit{Ars Conjectandi}
into Russian and the author of \citet{uspensky} brought the rigorous
probabilistic Russian tradition to the English speaking world after
moving to the United States.

Yakov Viktorovich Uspensky (1883--1947) is described by Markov in his
May, 1913, Foreward to the translation as ``Privat-Docent'' (roughly,
Assistant Professor) of St. Petersburg University. His academic contact
with Markov seems to have been through Markov's other great field of
interest, number theory. Uspensky's magisterial degree at this
university was conferred in 1910. He wrote on quadratic forms and
analytical methods in the additive theory of numbers. He was
``Privat-Docent'' 1912--1915, and Professor 1915--1923, and taught the
to-be-famous Russian number theorist I.M. Vinogradov in the Petrograd
incarnation of St. Petersburg. For his election to the Russian Academy
of Sciences in 1921, he had been nominated by A.A. Markov, V.A.
Steklov, and A.N. Krylov, and upto the time of elections in 1929 he was
the only mathematician in the Academy (\citet{bernoulli86}, p. 73). After
Markov's death in 1922, it was Uspensky who wrote a precis of Markov's
academic activity in the Academy's \textit{Izvestiia}, \textbf{17} (1923)
19--34. According to \citet{royden}, p. 243, the ``year 1929--1930 saw the
appointment of James Victor Uspenskly as an acting professor of
mathematics'' at Stanford University. He was professor of mathematics
there from 1931 until his death. He appears to have anglicized his name
and patronymic, Yakov Viktorovich, into James Victor, and it is under
this name, or just as J.V. Uspensky, that he appears in his
English-language writings. \citet{royden} writes that Uspensky had made
a trip to the U.S. in the early 1920s. When he did decide to come
permanently he came ``in style on a Soviet ship with his passage paid
for by the [Soviet] government'', which presumably was unaware of his intentions.

\section{Extensions. Necessary and sufficient conditions}

The expression (\ref{wlln1}) is the classical form of what is now
called the WLLN. We have confined ourselves to sufficient conditions
for (\ref{wlln1})
where $S_n =\sum_{i=1}^n X_i$ and the $ \{X_i, i=1,2, \ldots\}$ are
independent and not necessarily identically distributed. In particular,
in the tradition of Jacob Bernoulli's Theorem as limit theorem, we have
focused on the case of ``Bernoulli'' summands where $P(X_i=1) = p_i = 1
- P (X_i= 0)$.

Because of limitation of space we do not discuss the SLLN, and direct
the reader to our historical account (\citet{seneta92}), which begins with
\citet{borel} and \citet{cantelli17}. Further historical aspects may be
found in \citet{fisz63}, \citet{chung}, \citet{petrov} and \citet{krengel}.
The SLLN under ``Chebyshev's conditions'': $\{X_k \}$, $k= 1,2, \ldots$
pairwise independent, with variances well-defined and uniformly
bounded, was already available in \citet{rajchman}, but this source may
have been inaccessible to both Bernstein and Uspensky, not only because
of its language.

From the 1920s attention turned to \textit{necessary and sufficient}
conditions for the WLLN for independent summands. Kolmogorov in 1928
obtained the first such condition for ``triangular arrays'', and there
were a generalizations by Feller in 1937 and Gnedenko in 1944 (see
\citet{gnedkolm}, Section 22).

In Khinchin's paper on the WLLN in Cantelli's journal (\citet{khinchin36}) attention turns to necessary and sufficient conditions for the
existence of a sequence $\{ d_n\}$ of positive numbers such that
%
\begin{equation}
\label{wlln2} \frac{S_n}{d_n} \stackrel{p} {\to} 1 \mbox{ as } n \to\infty
\end{equation}
where the (i.i.d.) summands $X_i$ are
\textit{non-negative}.

Two new features in the consideration of limit theory for i.i.d.
summands make their appearance in Khinchin's several papers in
Cantelli's journal: a focus on the \textit{asymptotic structure of the
tails of the distribution function}, and the expression of this
structure in terms of what was later realized to be \textit{regularly
varying functions} (\citet{feller2}, \citet{seneta76}).

Putting $F(x) = P(X_i \leq x)$ and $ \nu(x) = \int_0^x (1-F(u))\,du$,
Khinchin's necessary and sufficient condition for (\ref{wlln2}) is
%
\begin{equation}
\label{wlln3} \frac{x(1 - F(x)}{\nu(x)} \to0 \mbox{ as } x \to\infty.
\end{equation}
This is equivalent to $\nu(x)$ being a slowly varying function at
infinity. In the event, $d_n$ can be taken as the unique solution of
$n\nu(d_n) =d_n$. A detailed account is given by
\citet{csorgosim}, Section 0. It is worth noting additionally that
$ \lim_{ x \to \infty} \nu(x) = E X_i \leq\infty$; and that if $x(1-
F(x)) = L(x)$, a slowly varying function at infinity (this includes the
case of finite mean $\mu$ when $\nu(x) \to\mu$), then (\ref{wlln3}) is
satisfied. (\citet{feller2} Section VII.7, p. 233, Theorem 3;
\citet{seneta76}, p. 87).

\citet{csorgosim}, motivated partly by the St. Petersburg
problem\footnote{Associated with Nicolaus and Daniel Bernoulli.}
consider, more generally than the WLLN sum structure, arbitrary linear
combinations of i.i.d. nonnegative random variables. When specializing
to sums $S_n$, however, they show in their Corollary 5 that
%
\begin{equation}
\label{wlln4} \frac{S_n}{n\nu(n)} \stackrel{p} {\to} 1
\end{equation}
if and only if
%
\begin{equation}
\label{wlln5} \frac{\nu(x \nu(x))}{\nu(x)} \to1 \mbox{ as } x \to\infty.
\end{equation}
They call (\ref{wlln5}) the \textit{Bojani\'c--Seneta condition.}\footnote
{It originates from \citet{bojanicsen}.}

Khinchin's Theorem itself was generalized by Feller (see, for example,
\citet{feller2} Section VII.7) in the spirit of \citet{khinchin36} for
i.i.d but not necessarily nonnegative summands.

\citet{petrov}, Chapter 6, Theorem 4, gives necessary and sufficient
conditions for the existence of a sequence of constants $\{b_n \}$ such
that $S_n/a_n-b_n \to0$ for any given sequence of positive constants
$\{a_n \}$ such that $ a_n \to\infty$, where the independent summands
$X_i$ are not necessarily identically distributed.

To conclude this very brief sketch, we draw the reader's attention to a
little-known necessary and sufficient condition for (\ref{wlln1}) to
hold, for arbitrarily dependent not necessarily identically distributed
random variables (see, for example, \citet{gnedenko63}).

\section*{Acknowledgements}

I am indebted to Bernard Bru and Steve Stigler for information, to an
anonymous Russian probabilist for suggesting emphasis on ``Markov's
Theorem'' and the contributions of S.N. Bernstein, and to the editors
Richard Davis and Thomas Mikosch for careful reading.



\printhistory

\end{document}